\documentclass[11pt]{article}
\usepackage{amsfonts}
\usepackage{eurosym}
\usepackage{latexsym}
\usepackage{amssymb}
\usepackage{graphicx}
\usepackage{amsmath}
\usepackage{float}
\usepackage{xcolor}

\newtheorem{theorem}{Theorem}

\newtheorem{corollary}[theorem]{Corollary}

\newtheorem{definition}[theorem]{Definition}

\newtheorem{lemma}[theorem]{Lemma}

\newtheorem{remark}[theorem]{Remark}

\floatstyle{ruled} \restylefloat{figure} \restylefloat{table}

\newcommand{\I}{\mathcal I}

\newcommand{\kom}[1]{}

\numberwithin{equation}{section}
\renewcommand{\kom}[1]{{\bf [#1]}}
\definecolor{darkjunglegreen}{rgb}{0.1, 0.14, 0.13}
\definecolor{blau}{rgb}{0.1,0.0,0.9}

\newcounter{komcounter}
\numberwithin{komcounter}{section}

\input{tcilatex}
\begin{document}

\title{Existence, uniqueness and regularity of solutions to systems of nonlocal obstacle problems related to optimal switching}
\author{%
Niklas L.P. Lundstr\"{o}m\thanks{Department of Mathematics and Mathematical Statistics, Ume\aa\ University,
SE-90187 Ume\aa , Sweden. E-mail: niklas.lundstrom@umu.se} \and
Marcus Olofsson\thanks{Department of Mathematics, Uppsala University, SE-75106 Uppsala, Sweden. E-mail: marcus.olofsson@math.uu.se} \and
Thomas \"{O}nskog\thanks{Department of Mathematics, KTH Royal Institute of Technology, SE-10044
Stockholm, Sweden. E-mail: onskog@kth.se}}
\maketitle

\begin{abstract}
\noindent
We study viscosity solutions to a system of nonlinear degenerate parabolic partial integro-differential equations with interconnected obstacles. This type of problem occurs in the context of optimal switching problems when the dynamics of the underlying state variable is described by an $n$-dimensional L\'{e}vy process. We first establish a continuous dependence estimate for viscosity sub- and supersolutions to the system under mild regularity, growth and structural assumptions on the partial integro-differential operator and on the obstacles and terminal conditions. Using the continuous dependence estimate, we obtain the comparison principle and uniqueness of viscosity solutions as well as Lipschitz regularity in the spatial variables. Our main contribution is construction of suitable families of viscosity sub- and supersolutions which we use as \textquotedblleft barrier functions" to prove H\"{o}lder continuity in the time variable, and, through Perron's method, existence of a unique viscosity solution. This paper generalizes parts of the results of Biswas, Jakobsen and Karlsen (2010) \cite{BJK10} and of Lundstr\"om, Nystr\"om and Olofsson (2014) \cite{LNO14, LNO14b} to hold for more general systems of equations.
\medskip


\medskip

\end{abstract}

\setcounter{equation}{0} \setcounter{theorem}{0}

\section{Introduction}

This paper deals with systems of partial integro-differential equations with
interconnected obstacles. The problem considered can be stated as
\begin{eqnarray}
\min \left\{ -\partial _{t}u_{i}\left( x,t\right) +\mathcal{F}_{i}\left(
x,t,u_{i}\left( x,t\right),Du_{i}\left( x,t\right) ,D^{2}u_{i}\left( x,t\right) ,u_{i}\left( \cdot
,t\right) \right) ,\frac{{}}{{}}\right. &&  \notag \\
\left. u_{i}\left( x,t\right) -\max\limits_{j\neq i}\left\{ u_{j}\left(
x,t\right) -c_{ij}\left( x,t\right) \right\} \right\} &=&0,
\label{problem} \\
u_{i}(x,T) &=&g_{i}\left( x\right) ,  \notag
\end{eqnarray}
for $x\in \mathbb{R}^{n}$, $t\in \lbrack 0,T]$ and $i\in \{1,2,...,m\}$.
Here the $\mathcal{F}_{i}$ are partial integro-differential operators of the
form
\begin{equation*}
\mathcal{F}_{i}\left( x,t,r,p,X,\varphi \left( \cdot \right) \right) =-%
\mathcal{L}_{i}\left( x,t,r,p,X\right) -\mathcal{J}_i\left( x,t,p,\varphi \left(
\cdot \right) \right) -f^{i}\left( x,t\right) ,
\end{equation*}
for local second order, possibly degenerate, linear operators $\mathcal{L}_{i}$, nonlocal operators $\mathcal{J}_i$ and functions $f^{i}$.

A particular case of problem \eqref{problem} arises in multi-modes switching
problems, in which the operators $\mathcal{F}_{i}$ and $\mathcal{I}_i$ can be interpreted as the
infinitesimal generators of the stochastic processes underlying the
optimization problem and where $c_{ij}$ is the cost of switching from state $%
i$ to state $j$. The viscosity solution of \eqref{problem} is then the value
function of the multi-modes switching problem from which one can find the
sought optimal strategy.

Since the pioneering work of Brennan and Schwarz \cite{BS85} dealing with a
two-modes switching problem describing the life cycle of an investment in
the natural resource industry, systems of variational inequalities with
interconnected obstacles of type \eqref{problem} has been extensively
studied. We here give a non-exhaustive list of contributions. Starting in the
purely local setting, i.e., when the operators $\mathcal{I}_{i} \equiv 0$,
problem \eqref{problem} and its connection to multi-modes switching have
been studied by, e.g., El-Asri and Hamad\`{e}ne \cite{AH09}, Djehiche, Hamad\`{e}ne, and
Popier \cite{DHP10}, Hamad\`{e}ne and Morlais \cite{HM12}, Hu and Tang \cite%
{HT07}, Lundstr\"{o}m, Nystr{\"{o}}m, and Olofsson \cite{LNO14b} and
Djehiche, Hamad\`{e}ne, Morlais and Zhao \cite{DHMZ14}.
In the majority of the above cited papers,
the main focus lies on existence and uniqueness of solutions to \eqref{problem}.

The nonlocal setting has not been studied to the same extent, but the amount
of literature is steadily growing. Biswas, Jakobsen, and Karlsen \cite{BJK10}
show that the dynamic programming principle for multi-modes switching
problems holds also for a nonlocal underlying stochastic process, i.e., when
the process is allowed to jump, and this is used to prove that the value
function of such a multi-modes switching problem indeed satisfies a system
like \eqref{problem}. The authors also proceed to prove existence,
uniqueness and some regularity of the viscosity solution. The nonlocal
setting has also been studied by Lundstr\"{o}m, Nystr{\"{o}}m, and Olofsson
\cite{LNO14}, using PDE techniques,
and recently by 
Hamad\`{e}ne and Zhao
\cite{HZ15}, Hamad\`{e}ne and Morlais \cite{HM16}
and Klimsiak \cite{K16, K16b} using mainly stochastic techniques.
In connection, important contributions to the mathematical theory of viscosity solutions to
second order partial integro-differential equations,
such as those considered in this paper,
was given by Barles and Imbert \cite{BI08} and Jakobsen and Karlsen \cite{JK06}.

A special technical difficulty when studying systems like \eqref{problem} is
the treatment of the obstacles. Therefore, many of the above references
choose to impose rather strong assumptions on the switching costs $c_{ij}$,
e.g., positive and constant switching costs. In fact, to the authors
knowledge, possibly negative switching costs are only treated by Asri and Fakhouri
\cite{AF12} and Martyr \cite{M14} using stochastic techniques and recently in \cite{LNO14} and
\cite{LNO14b} using PDE techniques. We note that it is natural to allow for
negative switching costs as these allow one to model the situation when, e.g., a government through environmental policies provides subsidies,
grants or other financial support to energy production facilities in case
they switch to more `green' production. In this case it is not a cost for the facility to switch, it is
a gain.

In this paper, we first establish a continuous dependence estimate (Theorem \ref%
{thm:continuousdependence}) by using standard methods for viscosity
solutions of integro PDEs. Such an estimate bounds the difference between a
solution to \eqref{problem} and a solution to a slightly different version
of \eqref{problem} in terms of the difference between the equations and the
terminal conditions. Continuous dependence estimates of this type are
important in themselves as they quantify the stability properties of
viscosity solutions, but are also useful in numerical analysis, see, e.g., Krylov
\cite{34BJK} and Barles and Jakobsen \cite{10BJK}. In this paper, the continuous dependence estimate is used to obtain
the comparison principle and uniqueness of solutions (Corollary \ref{cor:comparison}),
as well as regularity of the solution, both in space and time (Corollary \ref{cor:spatialreg} and Theorem \ref{thm:time}).
The main contributions of this paper are H\"{o}lder continuity in the time variable,
and the existence of a unique viscosity solution to system
\eqref{problem} (Theorem \ref{thm:existence}).
The proofs of these theorems rely on nontrivial constructions of
families of viscosity sub- and supersolutions (Lemma \ref{lemma:supersol} and Corollary \ref{lemma:subsol}) which we use as barrier functions in order to trap the viscosity solution via the comparison principle.

Our results generalize regularity and existence results of Biswas, Jakobsen, and
Karlsen \cite{BJK10} and of Lundstr\"{o}m, Nystr{\"{o}}m, and Olofsson \cite%
{LNO14} in the sense that we allow for more general systems of equations. We impose weaker assumptions on the operator and
less regularity and structural assumptions on the spatially dependent switching costs $c_{ij}$.
Our barrier functions also imply a more general existence result in the setting of Kolmogorov operators studied in Lundstr\"{o}m, Nystr{\"{o}}m, and Olofsson \cite{LNO14b}.


\setcounter{equation}{0} \setcounter{theorem}{0}

\section{Preliminaries}
\label{sec:assumptions}

In this section we introduce the notation used in the paper, discuss some
preliminaries and state the assumptions imposed on the system \eqref{problem}.\\

\noindent
For smooth functions $\varphi :\mathbb{R}^{n}\rightarrow \mathbb{R}$ we let
$D\varphi =(\partial _{x_{1}}\varphi ,\dots ,\partial _{x_{n}}\varphi )$
denote the spatial gradient of $\varphi $ and $D^{2}\varphi $ the Hessian
matrix of $\varphi$.
We denote the set of natural numbers by $\mathcal{N} = \{1,2,3,\dots\}$ and let $\mathcal{I}_{m}$ be the integer set $\{1,2,...,m\}$.
For any positive integer $p$ we  let $LSC_{p}\left( \mathbb{R}^{n}\times \lbrack0,T]\right)$ and $USC_{p}\left( \mathbb{R}^{n}\times \lbrack 0,T]\right)$
denote the spaces of lower- and upper semicontinuous functions,
respectively, on $\mathbb{R}^{n}\times \lbrack 0,T]$,
whose elements $h$ satisfy the growth condition
\begin{equation}\label{eq:growth-poly}
\left\vert h\left( x,t\right) \right\vert \leq K\left( 1+\left\vert
x\right\vert ^{p}\right),
\quad \text{whenever $\left( x,t\right) \in \mathbb{R}^{n}\times \lbrack 0,T]$.}
\end{equation}
Here and in the following, $K$ denotes a generic constant, $1 < K < \infty$, which may change value from line to line.
The space $C_p^{a,b}\left( \mathbb{R}^{n}\times \lbrack 0,T]\right)$ contains all functions
$\varphi :\mathbb{R}^{n}\times \lbrack 0,T]\rightarrow \mathbb{R}$ which are $a$ times continuously differentiable
in the spatial variables and $b$ times continuously differentiable in the time variable
and which satisfy the polynomial growth condition \eqref{eq:growth-poly}. 
From here on in, we fix the growth parameter $p$.
We denote the indicator function for the
closed unit ball in $\mathbb{R}^{l}$ by $\chi _{\{|z|\leq 1\}}$
and let $B(x,r)$ be the closed ball in $\mathbb{R}^n$ which has radius $r$ and is centered at $x$. 
We let $I_n$ denote the $n \times n$ identity matrix, 
and let $\mathbb{S}^{n}$ denote the space of $n\times n$ real symmetric matrices equipped with the positive semi-definite
ordering, i.e., for $X,Y\in \mathbb{S}^{n}$, we write $X\leq Y$ if $%
\langle (X-Y)\xi ,\xi \rangle \leq 0$ for all $\xi \in \mathbb{R}^{n}$.
 We will also make use of the matrix norm notation
\begin{equation*}
\left\Vert A\right\Vert :=\sup \{|\lambda |:\lambda \text{ is an eigenvalue of }A\}=\sup \{|\langle A\xi ,\xi \rangle |:|\xi |\leq 1,\xi \in \mathbb{R}^{n}\}.
\end{equation*}
The supremum norm is denoted
$\Vert h \Vert_{\infty}: = \sup_{x\in \mathbb{R}^n} |h(x)|$
for any function $h$ defined on $\mathbb{R}^n$.\\


We assume that the operator $\mathcal{F}_{i}$ can be decomposed into a local
second order operator $\mathcal{L}_{i}$,
a nonlocal integral operator $\mathcal{J}_i$ and a function $f^{i}$ such that
\begin{equation}
\mathcal{F}_{i}\left( x,t,r,p,X,\varphi \left( \cdot \right) \right)
=- \mathcal{L}_{i}\left( x,t,r,p,X\right)
 - \mathcal{J}_i\left( x,t,p,\varphi \left(\cdot \right) \right)
 - f^{i}\left( x,t\right),  \tag{$F_{1}$}  \label{Fass}
\end{equation}
for $\left( x,t,r,p,X\right) \in \mathbb{R}^{n}\times [0,T]\times\mathbb{R}\times \mathbb{%
R}^{n}\times \mathbb{S}^{n}$ and any smooth function $\varphi :$ $\mathbb{R}%
^{n}\rightarrow \mathbb{R}$. We assume that the local operators $\mathcal{L}_{i}$ can be
written as
\begin{equation*}
\mathcal{L}_{i}\left( x,t,r,p,X\right)
=\sum_{k,l=1}^{n} a^i_{kl}(x,t) X_{kl}  +  \sum_{k=1}^{n} b^i_{k}(x,t) p_{k} - c^i(x,t) r,
\end{equation*}
for continuous functions $a^i_{kl}$, $b^i_{k}$ and $c^i$.
We denote by $a^i$ the $n\times n$ matrix with elements $a^i_{kl}$,
and by $b^i$ the vector of length $n$ with elements $b^i_{k}$.
Moreover, we assume that
\begin{equation*}
a^i_{kl}(x,t) = (\sigma^i(x,t)(\sigma^i)^{\ast}(x,t))_{kl},
\quad \text{for } i \in \mathcal{I}_{m}, \quad k,l\in \mathcal{I}_{n},
\end{equation*}
for an $n\times n$ matrix $\sigma^i$ and where $(\sigma^i)^{\ast}$ is the transpose of $\sigma^i$.
The functions $\sigma^i_{kl}$, $b^i_{k}$, $c^i$ and $f^{i}$ are assumed to satisfy
\begin{align} \label{assump1elel}
&|\sigma^i_{kl}(x,t)-\sigma^i_{kl}(y,t)| + |b^i_{k}(x,t) - b^i_{k}(y,t)| + |c^i(x,t) - c^i(y,t)| \leq K|x-y|,\notag\\
&|f^i(x,t) - f^i(y,t)| \leq K \left( 1 + |x|^{p-1} + |y|^{p-1}\right) |x - y|\; \tag{$F_{2}$} \\
&|b^i_{k}(0,t)| + |\sigma^i_{kl}(0,t)| + |f^i(0,t)| - c^i(x,t) \leq K 
\notag
\end{align}
whenever $k,l\in \mathcal{I}_{n}$, $i \in \mathcal{I}_{m}$, $x,y\in \mathbb{R}^{n}$ and $t\in\lbrack 0,T]$. 
Concerning the nonlocal operators $\mathcal{J}_i$ we assume they can be written as
\begin{eqnarray*}
\mathcal{J}_i(x,t,p,\varphi \left( \cdot \right) )
&=& \int_{\mathbb{R}^{l}\setminus \{0\}}  \varphi \left( x+\eta^i \left(x,t,z\right) \right) - \varphi \left( x\right) -\chi_{\{|z|\leq1\}} \langle \eta^i\left( x,t,z\right), p \rangle \nu^i \left(dz\right),
\end{eqnarray*}
where $\nu^i$ is a positive Radon measure defined on $\mathbb{R}^{l}\setminus\{0\}$
and $\eta^i$ is an $\mathbb{R}^{n}$-valued function, continuous in $x$ and $t$ and Borel measurable in $z$.
We assume that $\nu^i$ and $\eta^i$ satisfy
\begin{align}
&\int_{0<|z|\leq 1}|z|^{2}\nu^i \left( dz\right) + \int_{|z|>1}|z|^{p}\nu^i \left( dz\right) < K ,  \tag{$F_{3}$}  \label{assump1elelnloc} \\
&|\eta^i_{k}\left( x,t,z\right) - \eta^i_{k}\left( y,t,z\right) | \leq K\left(1 + |z|\right)|x-y| \quad \textrm{and}\quad
|\eta^i_{k} \left( x,t,z\right) | \leq K \left(1 + |x|\right) |z|
\notag
\end{align}
whenever $k\in \mathcal{I}_{n}$, $i \in \mathcal{I}_{m}$, $x,y\in \mathbb{R}^{n}$, $t\in \lbrack 0,T]$ and $z\in \mathbb{R}^{l}$.

The functions $c_{ij}$ appearing in the obstacle are called \textquotedblleft switching costs" due to the connection between \eqref{problem} and optimal switching problems.
In light of this connection, the following definition makes sense.
\begin{definition}
A \textit{switching chain} from state $i$ to state $j$ is a sequence of
indices $(i_{1},\dots ,i_{l})\in \mathcal{I}_{m}^{l}$ such that $i_{1}=i$
and $i_{l}=j$. The set of switching chains from $i$ to $j$ is denoted $%
\mathcal{A}_{ij}$.
\end{definition}
We assume $c_{ij}$ to be continuous functions satisfying the classical no-loop condition, i.e.,
\begin{equation} \label{ass:noloop}
\min_{(i_{1},\dots ,i_{l})\in \mathcal{A}_{ii}}%
\sum_{k=1}^{l-1}c_{i_{k}i_{k+1}}(x,t)>0,  \tag{$O_{1}$}
\end{equation}
for all $i\in \mathcal{I}_{m}$ and $(x,t)\in \mathbb{R}^{n}\times \lbrack 0,T]$.
Moreover, we will need the stronger structural assumption
\begin{equation} \label{ass:structur}
c_{ik}(x,t) \leq c_{ij}(x,t)  + c _{jk}(x,t), 
\tag{$O_{2}$}
\end{equation}
whenever $i,j,k \in \I_m$, $x\in \mathbb{R}^{n}$ and $t \in [0,T]$. 
The assumption \eqref{ass:structur} is needed for our existence and time-regularity results,
(in particular, to prove Lemma \ref{lemma:supersol}),
but we stress that this assumption can be made without loss of generality in the context of optimal switching,
see Remark \ref{remark:swcosts}.

For Lemma \ref{lemma:supersol} we also need to assume that $c_{ij}$ is locally semi-concave in space,
locally Lipschitz continuous in both space and time, and satisfy a polynomial growth condition in space.
In particular, we assume that
\begin{align}\label{ass:reg}
&|c_{ij}(x,t)-c_{ij}(y,t')| \leq K \left( 1 + |x|^{p-1} + |y|^{p-1}\right) |(x,t) - (y,t')|, \notag\\
&D^{2}c_{ij}(z,s) \leq K \left( 1 + |z|^{p-2}\right) I_n, \tag{$O_{3}$}
\end{align}
whenever $i,j\in \mathcal{I}_{m}$, $x,y\in \mathbb{R}^{n}$, $t, t'\in \lbrack 0,T]$ and for almost every $(z,s)\in \mathbb{R}^{n}\times \lbrack 0,T]$.

Moreover, we assume that the terminal data $g_i$ is locally Lipschitz
continuous, and, to be able to achieve continuity up to the terminal time $T$,
that $g_{i}$ are consistent with the obstacle, i.e.,
\begin{align}\label{eq:gass2}
&|g_{i}(x)-g_{i}(y)|\leq K\left( 1 + |x|^{p-1} + |y|^{p-1}\right)|x-y|,
\quad g_{i}(x)\geq \max_{j\neq i}\left\{g_{j}\left( x\right) -c_{ij}(x,T)\right\},
\tag{$G$}  
\end{align}
whenever $i\in \mathcal{I}_{m}$ and $x,y\in \mathbb{R}^{n}$.

Since the matrices in the local operators $\mathcal{L}_{i}$ and the jump vectors $\eta^j$ are allowed to vanish,
we cannot expect any smoothing from the equation itself.
Therefore, a notion of weak solutions is needed and we will consider solutions in the viscosity sense.

\begin{definition}
\label{def:viscosity2} A vector $u=(u_{1},\dots ,u_{m})$, where $u_{i}\in
USC_{p}(\mathbb{R}^{n}\times \lbrack 0,T])$ (or $u_{i}\in LSC_{p}(\mathbb{R}%
^{n}\times \lbrack 0,T])$) for all $i\in \mathcal{I}_{m}$, is a viscosity
subsolution (supersolution) to system \eqref{problem} if $u_{i}(x,T)\leq g_{i}(x)$ (%
$u_{i}(x,T)\geq g_{i}(x)$) whenever $x\in \mathbb{R}^{n}$, $i\in \mathcal{I}%
_{m}$, and if the following holds. For every $(x_{0},t_{0})\in \mathbb{R}%
^{n}\times \lbrack 0,T)$ and $\varphi \in C_{p}^{2,1}(\mathbb{R}^{n}\times
\lbrack 0,T))$ such that $(x_{0},t_{0})$ is a global maximum (minimum) of $%
u_{i}-\varphi $, for some $i \in I_m$, we have
\begin{eqnarray*}
\min \left\{ -\partial _{t}\varphi_{i}\left( x_{0},t_{0}\right) +\mathcal{F}%
_{i}\left( x_{0}, t_{0}, u_{i}\left( x_{0},t_{0}\right), D\varphi \left( x_{0},t_{0}\right), D^{2}\varphi
\left( x_{0},t_{0}\right), \varphi \left( \cdot ,t_{0}\right)
\right) ,\frac{{}}{{}}\right. && \\
\left. u_{i}\left( x_{0},t_{0}\right) -\max\limits_{j\neq i}\left\{
u_{j}\left( x_{0},t_{0}\right) -c_{ij}\left( x_{0},t_{0}\right) \right\}
\right\} &\leq &(\geq )\,0,
\end{eqnarray*}
A vector $u=(u_{1},\dots ,u_{n})$, where $u_{i}\in C_{p}(\mathbb{R}^{n}\times \lbrack 0,T])$ for all $i\in \mathcal{I}_{m}$, is a viscosity solution to system \eqref{problem} if it is both a viscosity subsolution and a viscosity supersolution.
\end{definition}

Note that the test function appears in the nonlocal slot of the operator $\mathcal{I}_i$ in Definition \ref{def:viscosity2}. This is necessary due to the infinite activity of the jump measure $\nu^i$ close to the origin. However, away from the origin, the important property for $\mathcal{I}_i$ to be well-defined is not regularity but rather restrictions on its growth at infinity, see \eqref{assump1elelnloc}. Therefore, outside of the origin one may replace the test function $\varphi$ with the solution itself, $u_i$, and get an equivalent defintion of a viscosity solution, see Lemma 2.1 of \cite{BJK10}. When constructing barrier super- and subsolutions in Lemma \ref{lemma:supersol} and Corollary \ref{lemma:subsol}, we will use the above definition. However, when proving the continuous dependence estimate (Theorem \ref{thm:continuousdependence}) and the existence of a solution (Theorem \ref{thm:existence}), some necessary calculations follow those of \cite{BJK10} and \cite{LNO14b}. As \cite{BJK10} and \cite{LNO14b} use the latter defintion, we will do the same in the proofs of Theorem \ref{thm:continuousdependence} and Theorem \ref{thm:existence} to avoid repetition of lengthy calculations.


\setcounter{equation}{0} \setcounter{theorem}{0}

\section{Main results}
In this section, we list the main results of the paper. All proofs are
postponed to Section \ref{sec:proofs}. 
In the following, we write `depending on the data' to indicate dependence on (at most) the
constants $K$ and $p$ introduced in assumptions \eqref{Fass}--\eqref{assump1elelnloc}, %
\eqref{ass:noloop} --\eqref{ass:reg} and \eqref{eq:gass2}, as well as
dependence on the dimension $n$ and the terminal time $T$.
To state our first result,
which is a continuous dependence estimate,  we let,
for all $i \in \mathcal{I}_m$,
$\widehat{\mathcal{F}_i}$ denote the operator $\mathcal{F}_i$,
but with $\sigma^i_{kl}$, $b^i_{k}$, $c^i$, $f^{i}$, $\eta^i _{k}$ and $\nu^i$
replaced by $\widehat{\sigma}^i_{kl}$, $\widehat{b}^i_{k}$, $\widehat{c}^i$, $\widehat{f}^{i}$, $\widehat{\eta}^i_{k}$ and $\widehat{\nu}^i$,
respectively.

\begin{theorem}[Continuous dependence estimate]
\label{thm:continuousdependence}
Let $u=(u_{1},\dots ,u_{m})$ 
be a viscosity subsolution of system \eqref{problem} and let $\widehat{u}=(\widehat{u}_{1},\dots ,\widehat{u}_{n})$
be a viscosity supersolution of another system of the form \eqref{problem}
defined with $\widehat{\mathcal{F}}_{i}$, $\widehat{g}_{i}$ and $\widehat{c}_{ij}$
in place of $\mathcal{F}_{i}$, $g_{i}$ and $c_{ij}$.
Assume that both systems satisfy \eqref{Fass}--\eqref{assump1elelnloc}, \eqref{ass:noloop} and %
\eqref{eq:gass2}. Then there exists a positive constant $C$, depending only
on the data, such that
\begin{align*}
u_{i}(x,t)-\widehat{u}_{i}(x,t)
&\leq \, \max_{j\in \mathcal{I}_{m}} \Vert g_{j}-\widehat{g}_{j}\Vert_{\infty}
  + T \max_{j\in \mathcal{I}_{m}} \Vert f^{j}-\widehat{f}^{j}\Vert_{\infty}\\
&+ C  \max_{j\in \mathcal{I}_{m}}\bigg \{  \left\Vert c^j-\widehat{c}^j\right\Vert _{\infty}
+ \left\Vert b^j-\widehat{b}^j\right\Vert _{\infty}
  + \left\Vert \sigma^j-\widehat{\sigma}^j\right\Vert_{\infty}\\
&+ \left\Vert \int |\overline{\eta}^j|^{2}|\nu^j -\widehat{\nu}^j|(dz)\right\Vert_{\infty}^{1/2}
  + \left\Vert \int |\eta^j -\widehat{\eta}^j|^{2}\overline{\nu}^j(dz)\right\Vert_{\infty}^{1/2}\bigg
\},
\end{align*}
whenever $\left( x,t\right) \in \mathbb{R}^{n}\times \lbrack 0,T]$ and $i \in \mathcal{I}_m$, where $\overline{\eta}^i = \max \{\eta^i ,\widehat{\eta}^i\}$ and $\overline{\nu}^i =\max \{\nu^i ,\widehat{\nu}^i\}$.
\end{theorem}

\noindent
The classical comparison principle and Lipschitz regularity in the spatial
variables easily follows from Theorem \ref{thm:continuousdependence}. In
particular, setting $\widehat{\mathcal{F}}_{i}$ =$\mathcal{F}_{i}$,
$\widehat{g}_{i}=g_{i}$ and $\widehat{c}_{ij} = c_{ij}$ in Theorem \ref{thm:continuousdependence} gives the
following corollary.

\begin{corollary}[Comparison principle]
\label{cor:comparison}
Let $u^{-}=(u_{1}^{-},\dots ,u_{n}^{-})$ and $u^{+}=(u_{1}^{+},\dots ,u_{n}^{+})$ be a viscosity subsolution
and a viscosity supersolution of system \eqref{problem}, respectively.
Assume \eqref{Fass}--\eqref{assump1elelnloc}, \eqref{ass:noloop} and \eqref{eq:gass2}.
Then $u_{i}^{-}(x,t)\leq u_{i}^{+}(x,t)$ for all $\left( x,t\right) \in \mathbb{R}^{n}\times \lbrack 0,T]$
and $i\in \mathcal{I}_{m}$.
As a consequence, viscosity solutions to system \eqref{problem} are unique (in the class of polynomial growth).
\end{corollary}

\noindent
With the above results in place, consider a viscosity solution $u$ of system
\eqref{problem} and define,
for all $i,j\in \mathcal{I}_{m}$, $k,l\in
\mathcal{I}_{n}$ and $h\in \mathbb{R}^{n}$,
\begin{equation*}
\left( \widehat{\sigma}^i_{kl},\widehat{b}^i_{k},\widehat{c}^i,\widehat{f}^{i},\widehat{c}_{ij},\widehat{\eta}^i_{k}\right) \left( x,t\right)
:=\left( \sigma^i_{kl},b^i_{k},c^i,f^{i},c_{ij},\eta^i_{k}\right) \left( x+h,t\right).
\end{equation*}
Setting $\widehat{\nu}^i = \nu^i$ and $\widehat{g}_{i}\left( x\right) = g_{i}\left( x+h\right)$
it follows that $\widehat{u}\left( x,t\right)=u(x+h,t)$ is a viscosity solution to
\eqref{problem} with $\widehat{\mathcal{F}}_{i}$, $\widehat{g}_{i}$ and $\widehat{c}_{ij}$
in place of $\mathcal{F}_{i}$, $g_{i}$ and $\widehat{c}_{ij}$.
By the assumptions of Section \ref{sec:assumptions} it follows that we can bound the right-hand side of the
estimate in Theorem \ref{thm:continuousdependence} by $C\left(1 + |x|^{p-1} + |x + h|^{p-1}\right)|h|$,
where $C$ is a positive constant depending only on the data.
Hence, by an application of Theorem \ref{thm:continuousdependence} we have the following corollary.

\begin{corollary}[Lipschitz continuity in space]
\label{cor:spatialreg} Assume \eqref{Fass}--\eqref{assump1elelnloc}, %
\eqref{ass:noloop} and \eqref{eq:gass2}. Then there exists a constant $C$,
depending only on the data, such that for any viscosity solution $%
u=(u_{1},\dots ,u_{m})$ to system \eqref{problem} satisfying $u_{i}\in C_{p}\left(
\mathbb{R}^{n}\times \lbrack 0,T]\right) $ for all $i\in \mathcal{I}_{m}$,
it holds that
\begin{equation*}
|u_{i}(x,t)-u_{i}(y,t)|\leq C \left(1 + |x|^{p-1} + |y|^{p-1}\right) |x-y|,
\end{equation*}%
for any $i\in \mathcal{I}_{m}$, $x,y\in \mathbb{R}^{n}$ and $t\in \lbrack
0,T]$.
\end{corollary}

We proceed by stating our results on H\"{o}lder continuity in time and on existence of solutions. To prove these theorems we construct families of viscosity super- and subsolutions (Lemma \ref{lemma:supersol} and Corollary \ref{lemma:subsol}) which
we use as barrier functions in the comparison principle. To this end,
we need to impose the additional assumptions \eqref{ass:structur}--\eqref{ass:reg} on the
switching costs.

\begin{theorem}[H\"{o}lder continuity in time]
\label{thm:time} Assume \eqref{Fass}--\eqref{assump1elelnloc},
\eqref{ass:noloop}--\eqref{ass:reg} and \eqref{eq:gass2}. Then there
exists a constant $C$, depending only on the data, such that for any
viscosity solution $u=(u_{1},\dots ,u_{m})$ to system \eqref{problem} satisfying $%
u_{i}\in C_{p}\left( \mathbb{R}^{n}\times \lbrack 0,T]\right) $ for all $%
i\in \mathcal{I}_{m}$, it holds that
\begin{equation*}
|u_{i}(x,s)-u_{i}(x,t)|\leq C\left(1+|x|^{p}\right)|s-t|^{1/2}
\end{equation*}
for all $i\in \mathcal{I}_{m}$, $x\in \mathbb{R}^{n}$ and $s,t\in \lbrack 0,T]$.
\end{theorem}

\noindent
Finally, the following existence theorem is proved via Perron's method.
Here, the barrier functions from Lemma \ref{lemma:supersol} and Corollary \ref{lemma:subsol} are used to
ensure that the Perron solution is bounded and attains the terminal data.

\begin{theorem}[Existence]
\label{thm:existence} Assume \eqref{Fass}--\eqref{assump1elelnloc}, %
\eqref{ass:noloop}--\eqref{ass:reg} and \eqref{eq:gass2}. Then there
exists a unique viscosity solution $u=(u_{1},\dots ,u_{m})$ of system  %
\eqref{problem} satisfying $u_{i}\in C_{p}\left( \mathbb{R}^{n}\times
\lbrack 0,T]\right) $ for all $i\in \mathcal{I}_{m}$.
\end{theorem}

\begin{remark}
There is an $|x|^{p}$-dependence in the right hand side
of Theorem \ref{thm:time}, whereas the $|x|$-dependence in the corresponding
H\"{o}lder estimate in \cite{BJK10} (Lemma 5.3) is linear.
This is due to relaxed growth assumptions on $f_{i}$, $c_{ij}$ and $g_i$.
In particular, setting $p = 1$ in Theorem \ref{thm:time} we retrieve the result of
\cite{BJK10} in the more general setting studied here.
\end{remark}

\begin{remark}
Concerning generality we note that it should be possible to further relax the assumptions \eqref{Fass}--\eqref{assump1elelnloc}, by applying the full generality of the results of \cite{BI08} and \cite{JK06}. In particular, the continuous dependence estimate may be generalized using \cite{BI08} and \cite{JK06}. Given the validity of a continous dependence estimate, if the assumptions on the operator then implies Lipschitz continuity in space and the validity of \eqref{eq:final-bound},  then our barrier constructions hold and all our main results follows.
We have chosen to stay within the \textquotedblleft standard" assumptions \eqref{Fass}--\eqref{assump1elelnloc}
in this paper to avoid additional technicalities and  lenghty assumptions that are hard to interpret.
\end{remark}


\section{Proofs of the main results}
\label{sec:proofs}
In this section we prove Theorems \ref{thm:continuousdependence},
\ref{thm:time} and \ref{thm:existence}. \\

\noindent
{\bf Proof of Theorem \ref{thm:continuousdependence} (Continuous dependence estimate)}
We proceed along the lines of \cite[Theorem 5.1]{BJK10} to which we refer for additional details.

For constants
$\lambda, \theta, \gamma, \epsilon > 0$
we define the test function
\begin{align*}
\phi(t,x,y) = e^{\lambda(T - t)} \frac{\theta}{2} |x-y|^2 + e^{\lambda(T - t)} \frac{\epsilon}{2 + \gamma} \left( |x|^{2 + \gamma} + |y|^{2 + \gamma} \right),
\end{align*}
on $[0, T] \times \mathbb{R}^n \times \mathbb{R}^n$.
We double the variables by defining for $i \in \mathcal{I}_m$,
\begin{align*}
\Psi_i(t,x,y) = u_i(x, t) - \widehat u_i(y, t) - \phi(t,x,y) - \frac{\delta  (T-t)}{T} \sigma- \frac{\bar{\epsilon}}{t}
\end{align*}
where $0 < \delta < 1$, $\bar{\epsilon} > 0$, and
\begin{align*}
\sigma = \sup_{i,t,x,y} \bigg \{ u_i(x,t) - \widehat u_i(y,t) - \phi(t,x,y) - \frac{\bar{\epsilon}}{t}\bigg\} - \sigma_T,
\end{align*}
\begin{align*}
\sigma_T = \sup_{i,x,y} \bigg \{ u_i(x,T) - \widehat u_i(y,T) - \phi(T,x,y) - \frac{\bar{\epsilon}}{T}\bigg\}.
\end{align*}
From this we see that
\begin{align}\label{eq:cont-est_sigma_to_u}
u_i(x,t) - \widehat{u}_i(x,t) - e^{\lambda T} \epsilon |x|^{2 + \gamma} - \frac{\bar{\epsilon}}{t} \leq \sigma + \sigma_T, \quad \text{whenever}\quad \left( x,t\right) \in \mathbb{R}^{n}\times \lbrack 0,T], i\in \mathcal{I}_m,
\end{align}
and thus the main steps of the proof is to derive an upper bound on
$\sigma$ and $\sigma_T$. We start by establishing a bound for $\sigma$.
If $\sigma \leq 0$ we can take $0$ as the upper bound
and we are done. Therefore we will assume in the following that $\sigma > 0$.
By the upper semicontinuity of $u_i - \widehat{u}_i$,
the growth assumptions (provided $2 + \gamma > p$), and the penalization term $- \bar{\epsilon} / t$, there
exists $(i_0, t_0, x_0, y_0) \in \mathcal{I}_m \times (0, T] \times \mathbb{R}^n \times \mathbb{R}^n$ such that
\begin{align*}
\Psi_{i_0}(t_0,x_0,y_0) = \sup_{i,t,x,y} \Psi_{i}(t,x,y).
\end{align*}
The assumption $\sigma > 0$ forces $t_0 < T$.
To see this we observe that
\begin{align*}
\Psi_{i_0}(t_0,x_0,y_0) \geq \sup_{i,t,x,y}
\bigg\{u_i(x, t) - \widehat u_i(y, t) - \phi(t,x,y) - \frac{\bar{\epsilon}}{t} \bigg\} - \delta \sigma  = \sigma_T + (1-\delta) \sigma > \sigma_T,
\end{align*}
as $\delta < 1$,
while on the other hand $t_0 = T$ would imply $\Psi_{i_0}(t_0,x_0,y_0) = \sigma_T$.

Now we are in a position to apply the maximum principle for semicontinuous functions adapted to nonlocal systems.
As we allow switching costs $c_{ij}$ to depend on $(x,t)$, as well as polynomial growth of viscosity solutions, we may not apply \cite[Lemma 4.1]{BJK10} immediately.
However, we may apply the generalized version of this result found in \cite[Proof of Theorem 1.1]{LNO14b},
to retrieve the analogue of estimate $(5.3)$ in \cite{BJK10}. For each $0 < \kappa \leq 1$ there are symmetric matrices $X$ and $Y$,
and we can chose the index $i_0$ such that
\begin{align}\label{eq:cont-est_(5.3)}
-&\phi_t (t_0,x_0,y_0) + \frac{\delta \sigma}{T} + \frac{\bar{\epsilon}}{t_0^2}\notag \\
\leq &- \widehat{\mathcal{L}}_{i_0}(y_0, t_0, \widehat{u}^{i_0}(y_0,t_0), -D_y\phi(t_0,x_0,y_0), Y) - \widehat{f}^i(y_0,t_0) \notag\\
&-\widehat{\mathcal{J}}_{i_0,\kappa}(y_0, t_0, -D_y\phi(t_0,x_0,y_0), -\phi(t_0,x_0,\cdot))
-\widehat{\mathcal{J}}_{i_0}^{\kappa}(y_0, t_0, -D_y\phi(t_0,x_0,y_0), \widehat{u}^{i_0}(\cdot,t_0))\notag\\
& +\mathcal{L}_{i_0}(x_0, t_0, {u}^{i_0}(x_0,t_0), D_x\phi(t_0,x_0,y_0), X) +  {f}^i(x_0,t_0)\notag\\
&+\mathcal{J}_{i_0,\kappa}(x_0, t_0, D_x\phi(t_0,x_0,y_0), -\phi(t_0,\cdot,y_0))
+\mathcal{J}_{i_0}^{\kappa}(x_0, t_0, D_x\phi(t_0,x_0,y_0), {u}^{i_0}(\cdot,t_0)),
\end{align}
where the matrices $X$ and $Y$ satisfy standard upper bounds depending on the second derivatives of $\phi(t,x,y)$. 
We remind the reader that we here consider an alternative but equivalent definition of viscosity solutions and refer to 
\cite{BJK10} and \cite{LNO14b} for details. In \eqref{eq:cont-est_(5.3)}, the splitting of the nonlocal term is defined as in \cite[Definition 2.1]{BJK10}.

Following \cite{BJK10} we obtain the estimates
\begin{align}\label{eq:cont-est_1}
\text{tr}&\left( a^i(x_0,t_0) X\right) - \text{tr}\left( a^i(y_0,t_0) Y\right)\notag\\
&\leq C e^{\lambda(T - t_0)} \bigg\{ \theta |x_0 -y_0|^2 + \theta \Vert \sigma^i - \widehat{\sigma}^i\Vert ^2_{\infty}  + \epsilon \left( 1 + |x_0|^{2 + \gamma} + |y_0|^{2 + \gamma} \right)\bigg\},
\end{align}
\begin{align}\label{eq:cont-est_2}
\widehat{b}^i(y_0,t_0)& D_y \phi(t_0,x_0,y_0) + b^i(x_0,t_0) D_x \phi(t_0,x_0,y_0)\notag \\
&\leq C e^{\lambda(T - t_0)} \bigg\{ \theta |x_0 -y_0|^2 + \theta \Vert b^i - \widehat{b}^i\Vert ^2_{\infty}  + \epsilon \left( 1 + |x_0|^{2 + \gamma} + |y_0|^{2 + \gamma} \right)\bigg\},
\end{align}
\begin{align}\label{eq:cont-est_3}
|\widehat{f}^i(y_0,t_0) - {f}^i(x_0,t_0)| 
\leq  \Vert f^i - \widehat{f}^i\Vert _{\infty} + C \left(1 + |x_0|^{p-1} + |y_0|^{p-1} \right) |x_0 - y_0|,
\end{align}
for any $i \in \mathcal{I}_m$. Applying the polynomial growth assumptions of $u, \widehat u, c$ and $\widehat c$ yields
\begin{align}\label{eq:cont-est_4}
|\widehat{c}^i(y_0,t_0)& \widehat{u}(y_0,t_0) - {c}^i(x_0,t_0){u}(x_0,t_0)| \notag\\
&\leq C \left(1 + |x_0|^p + |y_0|^p\right) \Vert c^i - \widehat{c}^i\Vert _{\infty} + C \left(1 + |x_0|^{p} + |y_0|^{p} \right) |x_0 - y_0|.
\end{align}
For the nonlocal terms we obtain
\begin{align}\label{eq:cont-est_5}
&\mathcal{J}_{i_0,\kappa}(x_0, t_0, D_x\phi(t_0,x_0,y_0), -\phi(t_0,\cdot,y_0))
+\mathcal{J}_{i_0}^{\kappa}(x_0, t_0, D_x\phi(t_0,x_0,y_0), {u}^{i_0}(\cdot,t_0)),\notag\\
&-\widehat{\mathcal{J}}_{i_0,\kappa}(y_0, t_0, -D_y\phi(t_0,x_0,y_0), -\phi(t_0,x_0,\cdot))
-\widehat{\mathcal{J}}_{i_0}^{\kappa}(y_0, t_0, -D_y\phi(t_0,x_0,y_0), \widehat{u}^{i_0}(\cdot,t_0))\notag\\
& \leq
C \theta e^{\lambda(T - t_0)} \bigg \{ |x_0 - y_0|^2 + \left\Vert \int |\overline{\eta}^{i_0}|^{2}|\nu^{i_0} -\widehat{\nu}^{i_0}|(dz)\right\Vert_{\infty}
+ \left\Vert \int |\eta^{i_0} -\widehat{\eta}^{i_0}|^{2} \overline{\nu}^{i_0}(dz)\right\Vert_{\infty} \bigg \}\notag \\
&+ \mathcal{O}(\kappa) + C e^{\lambda(T - t_0)} \epsilon \left( 1 + |x_0|^{2 + \gamma} + |y_0|^{2 + \gamma} \right),
\end{align}
where $\overline{\eta}^{i_0} = \max \{\eta^{i_0}, \widehat{\eta}^{i_0}\}$ and $\overline{\nu}^{i_0} =\max \{\nu^{i_0} ,\widehat{\nu}^{i_0}\}$.
Now by \eqref{eq:cont-est_(5.3)}, estimates \eqref{eq:cont-est_1}-\eqref{eq:cont-est_5},
and the form of $\phi_t$, it follows that
\begin{align*}
&\lambda\left[ e^{\lambda(T - t_0)} \frac{\theta}2 |x_0 - y_0|^2 + e^{\lambda(T - t_0)} \frac{\epsilon}{2 + \gamma} \left(|x_0|^{2 + \gamma} + |y_0|^{2 + \gamma} \right) \right] + \frac{\delta \sigma}{T} + \frac{\bar{\epsilon}}{t_0^2}\notag\\
&\leq C \theta e^{\lambda(T - t_0)} \max_{i\in \mathcal{I}_m} \bigg \{ \Vert \sigma^i - \widehat{\sigma}^i\Vert ^2_{\infty} + \Vert b^i - \widehat{b}^i\Vert ^2_{\infty} + \left\Vert \int |\overline{\eta}^{i}|^{2}|\nu^{i} -\widehat{\nu}^{i}|(dz)\right\Vert_{\infty}
+ \left\Vert \int |\eta^{i} -\widehat{\eta}^{i}|^{2} \overline{\nu}^{i}(dz)\right\Vert_{\infty} \bigg \}\notag \\
& + \max_{i\in \mathcal{I}_m}  \Vert f^i - \widehat{f}^i\Vert _{\infty}
+ C \left(1 + |x_0|^p + |y_0|^p\right) \max_{i\in \mathcal{I}_m} \Vert c^i - \widehat{c}^i\Vert _{\infty}
+ C \theta e^{\lambda(T - t_0)}  |x_0 - y_0|^2\notag\\
& + C \left(1 + |x_0|^{p} + |y_0|^{p} \right) |x_0 - y_0|
+ C e^{\lambda(T - t_0)} \epsilon \left( 1 + |x_0|^{2 + \gamma} + |y_0|^{2 + \gamma} \right)
+ \mathcal{O}(\kappa)
\end{align*}
where the constant $C$ is not necessarily the same at each occurrence
but may depend only on the data.
In the above estimate, $t_0, x_0$ and $y_0$ are independent of $\kappa$,
so we can let $\kappa \to 0$ and ignore the term $\mathcal{O}(\kappa)$.
By taking $\lambda$ large enough, its magnitude depending only on the data,
we can conclude that
\begin{align}\label{eq:cont-est-long}
&\delta \sigma
\leq C T \theta \max_{i\in \mathcal{I}_m} \bigg \{ \Vert \sigma^i - \widehat{\sigma}^i\Vert ^2_{\infty} + \Vert b^i - \widehat{b}^i\Vert ^2_{\infty} + \left\Vert \int |\overline{\eta}^{i}|^{2}|\nu^{i} -\widehat{\nu}^{i}|(dz)\right\Vert_{\infty}
+ \left\Vert \int |\eta^{i} -\widehat{\eta}^{i}|^{2} \overline{\nu}^{i}(dz)\right\Vert_{\infty} \bigg \}\notag \\
& + T \max_{i\in \mathcal{I}_m} \Vert f^i - \widehat{f}^i\Vert _{\infty} + T \sup_{x,y}\Gamma(x,y),
\end{align}
where
\begin{align*}
\Gamma(x,y) &= C \left(1 + |x|^p + |y|^p\right) \max_{i\in \mathcal{I}_m} \Vert c^i - \widehat{c}^i\Vert _{\infty}
- \theta |x - y|^2
+ C \left(1 + |x|^{p} + |y|^{p} \right) |x - y|\notag\\
&- \epsilon \left( 1 + |x|^{2 + \gamma} + |y|^{2 + \gamma} \right)
+ C  \epsilon.
\end{align*}
In fact, by increasing $\lambda$ even further we see that we can take
\begin{align*}
\Gamma(x,y) &= C \left(1 + |x| + |y|\right)^p \max_{i\in \mathcal{I}_m} \Vert c^i - \widehat{c}^i\Vert _{\infty}
- \theta  |x - y|^2
+ C \left(1 + |x| + |y| \right)^{p} |x - y|\notag\\
&- \epsilon \left( 1 + |x| + |y| \right)^{2 + \gamma}
+ C \epsilon,
\end{align*}
and after a maximization with respect to $|x - y|$ we have
\begin{align*}
{\Gamma}(x,y) \leq  C \left(1 + |x| + |y|\right)^p \max_{i\in \mathcal{I}_m} \Vert c^i - \widehat{c}^i\Vert _{\infty}
+ C \frac{ \left(1 + |x| + |y| \right)^{2p}}{\theta } -  \epsilon \left( 1 + |x| + |y| \right)^{2 + \gamma}
+ C \epsilon,
\end{align*}
and so
\begin{align*}
{\Gamma}(x,y) \leq  C \left(1 + |x| + |y|\right)^{2p} \left( \max_{i\in \mathcal{I}_m} \Vert c^i - \widehat{c}^i\Vert _{\infty}
+ \frac{1}{\theta} \right) -  \epsilon \left( 1 + |x| + |y| \right)^{2 + \gamma}
+ C \epsilon.
\end{align*}
Now, pick $\gamma = 4p - 2$ and maximize anew, this time with respect to $\left( 1 + |x| + |y| \right)$.
The result is
\begin{align*}
{\Gamma}(x,y) \leq \frac{C}{\epsilon} \left( \max_{i\in \mathcal{I}_m} \Vert c^i - \widehat{c}^i\Vert ^2_{\infty}
+ \frac{1}{\theta^2} \right) + C \epsilon.
\end{align*}
and by now choosing $\epsilon = 1 / \theta$ we can conclude that
\begin{align}\label{eq:Gamma_bound}
{\Gamma}(x,y) \leq {C} \theta \max_{i\in \mathcal{I}_m} \Vert c^i - \widehat{c}^i\Vert ^2_{\infty} + \frac{C}{\theta}.
\end{align}

We next estimate $\sigma_T$.
We have, using \eqref{eq:gass2}, that
\begin{align*}
\sigma_T &= \sup_{i,x,y} \bigg \{ g_i(x) - \widehat g_i(y) - \phi(T,x,y) - \frac{\bar{\epsilon}}{T}\bigg\}\notag\\
&\leq \max_{i\in \mathcal{I}_m} \Vert g^i - \widehat{g}^i\Vert _{\infty}
+ \sup_{x,y} \bigg \{ C \left( 1 + |x|^{p-1} + |y|^{p-1}\right) |x-y| - \phi(T,x,y) \bigg\}\notag\\
&=  \max_{i\in \mathcal{I}_m} \Vert g^i - \widehat{g}^i\Vert _{\infty} + \sup_{x,y} \Gamma_T(x,y),
\end{align*}
where $\Gamma_T(x,y)$ can be bounded in a similar way as $\Gamma(x,y)$, i.e.,
\begin{align*}
\Gamma_T(x,y) &\leq C \left( 1 + |x| + |y|\right)^{p-1} |x-y| -  \frac{\theta}{2} |x-y|^2 - \frac{\epsilon}{C} \left( 1 + |x| + |y| \right)^{2 + \gamma}\notag\\
&\leq C \frac{ \left(1 + |x| + |y| \right)^{2p}}{\theta } - \frac{\epsilon}{C} \left( 1 + |x| + |y| \right)^{4p} \leq \frac{C}{\theta^2 \epsilon} = \frac{C}{\theta}
\end{align*}
and hence
\begin{align}\label{eq:sigma_T_bound}
\sigma_T \leq \max_{i\in \mathcal{I}_m} \Vert g^i - \widehat{g}^i\Vert _{\infty} + \frac{C}{\theta}.
\end{align}

Collecting the estimates \eqref{eq:cont-est-long}, \eqref{eq:Gamma_bound} and \eqref{eq:sigma_T_bound}, 
sending $\delta \to 1$ and inserting them in \eqref{eq:cont-est_sigma_to_u} yields,
after noting that the term $-e^{\lambda T} \epsilon |x|^{2 + \gamma}$ in \eqref{eq:cont-est_sigma_to_u} can be absorbed by $\Gamma(x,y)$
(by an increase in $\lambda$),
\begin{align*}
u_i(x,t) &- \widehat{u}_i(x,t) - \frac{\bar{\epsilon}}{t} \leq \sigma + \sigma_T + e^{\lambda T} \epsilon |x|^{2 + \gamma}\notag\\
&\leq C T \theta \max_{j\in \mathcal{I}_m} \bigg \{ \Vert \sigma^j - \widehat{\sigma}^j\Vert ^2_{\infty} + \Vert b^j - \widehat{b}^j\Vert ^2_{\infty} + \Vert c^j - \widehat{c}^j\Vert ^2_{\infty} + \left\Vert \int |\overline{\eta}^{j}|^{2}|\nu^{j} -\widehat{\nu}^{j}|(dz)\right\Vert_{\infty}\notag\\
&+ \left\Vert \int |\eta^{j} -\widehat{\eta}^{j}|^{2} \overline{\nu}^{j}(dz)\right\Vert_{\infty} \bigg \}
+ \max_{j\in \mathcal{I}_m} \Vert g^j - \widehat{g}^j\Vert _{\infty} + T \max_{j\in \mathcal{I}_m} \Vert f^j - \widehat{f}^j\Vert _{\infty} + \frac{C T}{\theta},
\end{align*}
whenever $\left(x,t\right) \in \mathbb{R}^{n}\times \lbrack 0,T]$
and $i\in \mathcal{I}_m$.
After minimizing the right hand side with respect to $\theta$ and sending $\bar{\epsilon} \to 0$, we have
\begin{align*}
&u_i(x,t) - \widehat{u}_i(x,t) \leq C T \max_{j\in \mathcal{I}_m} \bigg \{ \Vert \sigma^j - \widehat{\sigma}^j\Vert _{\infty} + \Vert b^j - \widehat{b}^j\Vert _{\infty} + \Vert c^j - \widehat{c}^j\Vert _{\infty} + \left\Vert \int |\overline{\eta}^{j}|^{2}|\nu^{j} -\widehat{\nu}^{j}|(dz)\right\Vert^{1/2}_{\infty}\notag\\
&+ \left\Vert \int |\eta^{j} -\widehat{\eta}^{j}|^{2} \overline{\nu}^{j}(dz)\right\Vert^{1/2}_{\infty} \bigg \}
+ \max_{j\in \mathcal{I}_m} \Vert g^j - \widehat{g}^j\Vert _{\infty} + T \max_{j\in \mathcal{I}_m} \Vert f^j - \widehat{f}^j\Vert _{\infty},
\end{align*}
whenever $\left( x,t\right) \in \mathbb{R}^{n}\times \lbrack 0,T]$ and $i\in \mathcal{I}_m$.
This completes the proof of the theorem. $\hfill \Box$\\


We will now prove Theorem \ref{thm:time} and Theorem \ref{thm:existence} by building appropriate barrier functions.
The barrier functions will be constructed as families of viscosity super- and subsolutions to \eqref{problem} which, by the comparison principle (Corollary \ref{cor:comparison}), will give bounds for the unique viscosity solution from above and below, respectively.
Before going into the proof we note that the difficulty lies in constructing an appropritate family of supersolutions which exceed the obstacle.
Our main idea for this construction is to include the switching costs explicitly in the barrier.
Since the switching costs are allowed to be non-smooth, the operator cannot be applied directly and we need to consider approximation arguments and viscosity solution theory.
A suitable family of subsolutions can be constructed independent of the obstacle, and is therefore much simpler.
Lemma \ref{lemma:supersol} gives an appropriate family of viscosity supersolutions.
The construction of this family was inspired by related arguments in
Lundstr\"om, Nystr\"om and Olofsson \cite{LNO14,LNO14b}, Biswas, Jacobsen and Karlsen \cite{BJK10},
Ishii and Sato \cite{IS04} and Lundstr\"om and \"Onskog \cite{L�1�7}.

\begin{lemma}[Upper barrier]
\label{lemma:supersol}
Assume \eqref{Fass}--\eqref{assump1elelnloc}, \eqref{ass:noloop}--\eqref{ass:reg} and let
$h=(h_{1},...,h_{m}):$ $\mathbb{R}^{n}\rightarrow \mathbb{R}^{m}$ be a function satisfying \eqref{eq:gass2} for some positive integer $p$.
Then there exists a positive constant $c$, depending only on $K,n$ and $p$,
such that for all $\left( y,s\right) \in \mathbb{R}^{n}\times \left[ 0,T\right]$
and all $i\in \mathcal{I}_{m}$,
the function $\psi^{i,y,s}=(\psi _{1}^{i,y,s},...,\psi _{m}^{i,y,s})$,
defined as
\begin{equation*}
\psi _{j}^{i,y,s}(x,t)
= c\, \lambda\, e^{c(s-t)} \left\{ A\,c \left(s - t\right) + \frac{A}{\lambda^2} + B |x-y|^{2} + |x-y|^{p}\right\}
+ h_{i}(y)
+ c_{ij}(x,t)
\end{equation*}
for all $j\in \mathcal{I}_{m}$, where $A = \left( 1 + |y|^{p}\right)$ and $B = \left( 1 + |y|^{p-2}\right)$,
is a viscosity supersolution to \eqref{problem} with terminal condition given by $h$, in $\mathbb{R}^{n}\times \lbrack 0,s)$ and whenever $\lambda \geq 1$.
\end{lemma}

\begin{corollary}[Lower barrier]
\label{lemma:subsol}
Assume \eqref{Fass}--\eqref{assump1elelnloc}, \eqref{ass:noloop}--\eqref{ass:reg} and let
$h, c, A$ and $B$ be as in Lemma \ref{lemma:supersol}.
Then for all $\left( y,s\right) \in \mathbb{R}^{n}\times \left[ 0,T\right]$
and all $i\in \mathcal{I}_{m}$,
the function $\check \psi^{y,s} = (\check \psi _{1}^{y,s},...,\check \psi _{m}^{y,s})$ defined as
\begin{equation*}
\check \psi_j^{y,s}(x,t)
= - c\, \lambda\, e^{c(s-t)} \left\{ A\,c \left(s - t\right) + \frac{A}{\lambda^2} + B |x-y|^{2} + |x-y|^{p}\right\} + h_{j}(y),
\end{equation*}
is a viscosity subsolution to \eqref{problem} with terminal condition given by $h$, in $\mathbb{R}^{n}\times \lbrack 0,s)$ and
whenever $\lambda \geq 1$.
\end{corollary}

\noindent
{\bf Proof.}
This follows by repeating steps 1 and 3 in the proof of Lemma \ref{lemma:supersol} given below.
Both steps are simpler in this case since $\check \psi_j^{y,s}$ does not involve the switching costs $c_{ij}$.   $\hfill\Box$\\

\noindent
Before proving Lemma \ref{lemma:supersol} we recall two well-known results in 
Lemma \ref{lemma:Jensen} and Lemma \ref{lemma:C11-semikonvexsemikoncave}, needed when we prove that our 
family of functions in Lemma \ref{lemma:supersol} consists of viscosity supersolutions.

\begin{lemma}\label{lemma:Jensen}
Let $\varphi : \mathbb{R}^n \to \mathbb{R}$ be semiconvex and $\widehat x$ be a strict local maximum point of $\varphi$.
For $p \in \mathbb{R}^n$, set $\varphi_p(x) = \varphi(x) + \langle p, x \rangle$.
Then for $r, \delta > 0$,
$$
K = \{ x \in B(\widehat x, r): \text{there exists $p \in B(0,\delta)$ for which $\varphi_p$ has a local maximum at $x$} \}
$$
has positive measure.
\end{lemma}

\noindent
{\bf Proof.}
This result is given as Lemma A.3 in \cite{CIL92} to which we refer for a proof. $\hfill\Box$

\begin{lemma}\label{lemma:C11-semikonvexsemikoncave}
Let $\Omega \subset \mathbb{R}^n$. 
A function $\varphi : \Omega \to \mathbb{R}$ is differentiable with derivative satisfying
$|D\varphi(x) - D\varphi(y)| \leq K |x-y|$ for all $x,y \in \Omega$
if and only if $\varphi$ is both $K$-semiconvex and $K$-semiconcave on $\Omega$,
i.e., both $\varphi(x) + \frac{K}{2}|x|^2$ and $-\varphi(x) + \frac{K}{2}|x|^2$ are convex functions on $\Omega$.
\end{lemma}

\noindent
{\bf Proof.}
A proof can be found in, e.g., Harvey and Lawson \cite[Theorem A.1]{HL13}. $\hfill\Box$\\

\noindent
Armed with Lemmas \ref{lemma:Jensen} and \ref{lemma:C11-semikonvexsemikoncave}
we are ready to prove Lemma \ref{lemma:supersol}.\\

\noindent
{\bf Proof of Lemma \ref{lemma:supersol} (Upper barrier)}
The proof naturally split into three steps.

\emph{Step 1: $\psi _{j}^{i,y,s}$ satisfies the terminal condition.}
We have to show that
\begin{align}\label{eq:initial}
\psi _{j}^{i,y,s}(x,s) \geq h_j(x) \quad \text{whenever} \quad x \in \mathbb{R}^n.
\end{align}
Using the fact that $ab\leq {a^{2}+b^{2}}$ we see that, for all $\lambda \geq 1$,
\begin{eqnarray*}
\psi _{j}^{i,y,s}(x,s)
&=&  c \left\{ \lambda\,B\,|x-y|^{2} + \frac{A}{\lambda}  +\lambda\, |x-y|^{p} \right\}
+ h_{i}(y)
+ c_{ij}(x,s)\\
&\geq& c \left\{  \sqrt{AB} |x-y| + \lambda\, |x-y|^{p} \right\}
+ h_{i}(y)
+ c_{ij}(x,s)\\
&\geq&  \frac{c}2 \left\{   1 + |y|^{p-1} +  |x-y|^{p-1} \right\} |x-y|
+ h_{i}(y)
+ c_{ij}(x,s).
\end{eqnarray*}
Now for any $k > 0$ it holds that
%
\begin{align}\label{eq:trivialellertekniskt}
 |y|^{k}  +  |x-y|^{k} \geq \frac{1}{2^{k+1}} \left( |x|^{k} + |y|^{k}\right)
\end{align}
and, using this inequality and that $h_{i}$ satisfies \eqref{eq:gass2} we can conclude that for $c \geq K 2^{p+1} $
\begin{eqnarray*}
\psi _{j}^{i,y,s}(x,s)
&\geq&  \frac{c}{2^{p+1}} \left( 1 + |x|^{p-1} + |y|^{p-1}\right) |x-y|
+ h_{i}(y)
+ c_{ij}(x,s)\\
&\geq& h_{i}(x)+c_{ij}(x,s) \geq h_{j}(x),  
\end{eqnarray*}
for all $\lambda \geq 1$, all $j\in \mathcal{I}_{m}$, and all $x\in \mathbb{R}^{n}$.
This proves \eqref{eq:initial} and therefore the terminal condition is fulfilled.\\

\emph{Step 2: $\psi _{j}^{i,y,s}$ exceeds the obstacle.}
To show that $\psi _{j}^{i,y,s}$ exceeds the obstacle we have to show that
\begin{equation}\label{eq:combine2}
\psi _{j}^{i,y,s}(x,t)-\max_{k\neq j}\left\{ \psi
_{k}^{i,y,s}(x,t)-c_{jk}(x,t)\right\} \geq 0
\end{equation}
at all points $(x,t)\in \mathbb{R}^{n}\times \lbrack 0,s]$ and whenever $i,j,k \in \mathcal{I}_m$.
To do so we note that assumption \eqref{ass:structur} reads
\begin{equation*}
c_{ij}(x,t)+c_{jk}(x,t)-c_{ik}(x,t)\geq 0,
\end{equation*}
for all $i,j,k \in \mathcal{I}_m$, and hence
\begin{equation*}
\psi _{j}^{i,y,s}(x,t)-(\psi
_{k}^{i,y,s}(x,t)-c_{jk}(x,t))=c_{ij}(x,t)-(c_{ik}(x,t)-c_{jk}(x,t))\geq 0.
\end{equation*}
%
Therefore, inequality \eqref{eq:combine2} is satisfied and our supersolution candidate exceeds the obstacle.\\

\emph{Step 3: $\psi _{j}^{i,y,s}$ satisfies the equation of being a supersolution.}
We cannot apply the operator
directly to $\psi _{j}^{i,y,s}$ since the switching costs $c_{ij}$ are in
general not differentiable.
Instead, we consider a viscosity solution
approach.
According to Definition \ref{def:viscosity2} we need to show that,
for $c$ large enough,
\begin{equation}\label{eq:eqpart}
-\partial _{t}\varphi\left( \widehat{x},\widehat{t}\right) +\mathcal{F}_{j}\left( \widehat{x},\widehat{t},\psi _{j}^{i,y,s}\left( \widehat{x},\widehat{t}\right),D\varphi\left( \widehat{x},\widehat{t}\right) ,D^{2}\varphi\left( \widehat{x},\widehat{t}\right) ,\varphi\left( \cdot ,\widehat{t}\right) \right) \geq \,0,
\end{equation}
whenever $\varphi $ is a $C_{p}^{2,1}(\mathbb{R}^{n}\times \lbrack 0,s))$
function such that,
for some $j \in I_m$,
$\psi _{j}^{i,y,s}-\varphi$ has a global minimum at $(\widehat{x},\widehat{t}) \in \mathbb{R}^n \times [0,s)$.
We may w.l.o.g. assume that the minimum is strict.

Using the notation
\begin{equation*}
\psi_{j}^{i,y,s}(x,t) = \gamma^{i,y,s}(x,t) + c_{ij}(x,t),
\end{equation*}
where
\begin{equation*}
\gamma^{i,y,s}(x,t) = c\, \lambda\, e^{c(s-t)} \left\{ A\,c \left(s - t\right) + \frac{A}{\lambda^2} + B |x-y|^{2} + |x-y|^{p}\right\} + h_{i}(y),
\end{equation*}
we find from assumption \eqref{ass:reg} that
\begin{align*}
\partial_t \varphi(\widehat{x},\widehat{t}) \leq \partial_t \gamma^{i,y,s}(\widehat{x},\widehat{t})
 + K\left(1 + |\widehat{x}|^p\right)
\quad \text{and} \quad
| D_{x_k} \varphi(\widehat{x},\widehat{t})| \leq D_{x_k}\gamma^{i,y,s}(\widehat{x},\widehat{t}) + K\left(1 + |\widehat{x}|^{p-1}\right),
\end{align*}
for all $k \in \mathcal{I}_n$.
Thus
\begin{align*}
-\partial_t \varphi(\widehat{x},\widehat{t})
&\geq c^2\,\lambda\, e^{c (s-\widehat{t})}  \left\{ A + B |\widehat{x}-y|^{2} + |\widehat{x}-y|^{p}\right\}
- K\left( 1 + |\widehat{x}|^p \right),
\end{align*}
and, as \eqref{eq:trivialellertekniskt} implies
\begin{align*}
\left\{ A + |\widehat{x}-y|^{p}\right\}
&\geq \frac{1}{2}\left\{ 1 + |y|^{p} + |y|^{p} + |\widehat{x}-y|^{p}\right\}
\geq \frac{1}{2^{p+2}} \left\{ 1 + |y|^{p} + |\widehat{x}|^{p}\right\},
\end{align*}
we conclude that
\begin{align}\label{eq:dt-bound}
-\partial_t \varphi(\widehat{x},\widehat{t})
&\geq c^2\,\lambda\, e^{c (s-\widehat{t})}  \frac{1}{2^{p+2}} \left\{ 1 + |\widehat{x}|^{p} + |y|^{p} \right\}
- K\left( 1 + |\widehat{x}|^p \right),
\end{align}
whenever $(\widehat{x},\widehat{t}) \in \mathbb{R}^n \times [0,T]$.
Similarly, for the first derivative in space we have that
\begin{align}\label{eq:first_der}
\left| D \varphi(\widehat{x},\widehat{t}) \right|
&\leq c\, \lambda\, e^{c(s-\widehat{t})} \left\{ 2 B |\widehat{x}-y| + p |\widehat{x}-y|^{p-1} \right\}
+ K \left( 1 + |\widehat{x}|^{p-1} \right) \notag\\
&\leq c\, \lambda\, e^{c(s-\widehat{t})} C \left( 1 + |\widehat{x}|^{p-1} + |y|^{p-1}\right)
\end{align}
since
\begin{align*}
\frac{\partial}{\partial x_k} \gamma^{i,y,s}(x,t)
= c\, \lambda\, e^{c(s-t)} \left\{ 2 B + p |x-y|^{p-2} \right\} \left(x_k - y_k\right).
\end{align*}
Here and in the following,
by $C$ we will denote a constant, $1\leq C<\infty$,
not necessarily the same at each occurrence,
which may depend only on $K, n$ and $p$.

To establish an upper bound for $D^2\varphi$, 
we first note that since $\varphi$ is twice differentiable in space,
it is, by Lemma \ref{lemma:C11-semikonvexsemikoncave}, locally semi-concave in space.
Hence, by assumption \eqref{ass:reg} the function $\psi _{j}^{i,y,s}-\varphi$ is also semi-concave.
We can thus apply Lemma \ref{lemma:Jensen} to obtain a sequence
$(x_{k}, q_{k}) \in \mathbb{R}^{n}\times \mathbb{R}^n$
such that $x_{k} \rightarrow \widehat{x}$,
$ q_{k}\rightarrow 0$ and such that the function
\begin{equation*}
\psi _{j}^{i,y,s}(x,\widehat{t}) - \varphi (x,\widehat{t})  + \langle q_{k}, x \rangle
\end{equation*}
has a local minimum at $(x_{k},\widehat{t})$.
Furthermore, from Lemma \ref{lemma:Jensen} it follows that we may assume that
$D^{2}\psi _{j}^{i,y,s}(x_{k},\widehat{t})$ exists for all $k$.
Therefore,
\begin{equation*}
D^{2}\varphi (x_{k},\widehat{t})
\leq D^{2}\psi_{j}^{i,y,s}(x_{k},\widehat{t})
= D^{2}\gamma^{i,y,s}(x_{k},\widehat{t}) + D^{2}c_{ij}(x_{k},\widehat{t}),
\end{equation*}
for all $i,j\in \mathcal{I}_{m}$.
Taking the limit as $k\rightarrow \infty$ and using \eqref{ass:reg} then yields
\begin{equation}\label{eq:matrixinequality_secondderivative}
D^{2}\varphi(\widehat{x},\widehat{t})
\leq D^{2}\gamma^{i,y,s}(\widehat{x},\widehat{t}) + K \left( 1 + |\widehat{x}|^{p-2}\right) I_{n}.
\end{equation}
%
%
Now, let $\delta_{kl}$ denote the Kronecker-delta and observe that
\begin{align*}
\frac{\partial^2}{\partial x_k \partial x_l} \gamma^{i,y,s}(x,t)
&= c\, \lambda\, e^{c(s-t)} \left\{ 2 B + p |x-y|^{p-2} \right\} \delta_{kl} \\
&+ c\, \lambda\, e^{c(s-t)} p\left(p-2\right) |x-y|^{p-4} \left(x_l - y_l\right) \left(x_k - y_k\right),
\end{align*}
from which we get, using \eqref{eq:matrixinequality_secondderivative}, that
\begin{align}\label{eq:secondderivativebound}
\Vert  D^{2}\varphi(\widehat{x},\widehat{t})\Vert
&\leq c\, \lambda\, e^{c(s-\widehat{t})} \left\{ 2 B + p(p-1) |\widehat{x}-y|^{p-2} \right\}
+ K \left( 1 + |\widehat{x}|^{p-2}\right) \notag \\
&\leq c\, \lambda\, e^{c(s-\widehat{t})} C \left(1 + |\widehat{x}|^{p-2} + |y|^{p-2} \right).
\end{align}

Next, by noting that \eqref{assump1elel} implies
$|a^j_{kl}(\widehat{x},\widehat{t})| \leq C \left(1 + |\widehat{x}|^2\right)$,
$c^j(\widehat{x},\widehat{t}) \geq - K$
and $|b^j_{k}(\widehat{x},\widehat{t})| \leq K \left(1 + |\widehat{x}|\right)$,
we conclude that
\begin{align}\label{eq:L-bound}
-&\mathcal{L}_{j}\left(\widehat{x},\widehat{t},\psi _{j}^{i,y,s}\left( \widehat{x},\widehat{t}\right),D\varphi(\widehat{x},\widehat{t}),D^2\varphi(\widehat{x},\widehat{t})\right)\notag\\
= & - \text{tr}\left( a^j(\widehat{x},\widehat{t}) D^2\varphi(\widehat{x},\widehat{t}) \right)
-  \sum_{k=1}^{n} b^j_{k}(\widehat{x},\widehat{t}) \left(D\varphi(\widehat{x},\widehat{t})\right)_{k}
+ c^j(\widehat{x},\widehat{t}) \psi _{j}^{i,y,s}(\widehat{x},\widehat{t})\notag\\
\geq & - n \Vert a^j(\widehat{x},\widehat{t})\Vert \, \Vert D^2\varphi(\widehat{x},\widehat{t})\Vert
- |b^j_{k}(\widehat{x},\widehat{t})|\, |D\varphi(\widehat{x},\widehat{t})|
- K |\psi _{j}^{i,y,s}(\widehat{x},\widehat{t})|\notag\\
\geq &  - C \left(1 + |\widehat{x}|^2\right)  c\, \lambda\, e^{c(s-\widehat{t})} C \left(1 + |\widehat{x}|^{p-2} + |y|^{p-2} \right)                            		- K \left(1 + |\widehat{x}|\right)    c\, \lambda\, e^{c(s-\widehat{t})} C \left( 1 + |\widehat{x}|^{p-1} + |y|^{p-1}\right)\notag\\
	&    - K c\, \lambda\, e^{c(s-\widehat{t})} C \left( 1 + |\widehat{x}|^{p} + |y|^{p}\right)\notag\\
\geq & -c\, \lambda\, e^{c(s-\widehat{t})} C \left( 1 + |\widehat{x}|^{p} + |y|^{p}\right).
\end{align}

Regarding the nonlocal term $\mathcal{J}_j$, we first note that since $\varphi$ is $C_p^{2,1}(\mathbb{R}^{n}\times \lbrack 0,s))$,
it follows from Taylor's theorem that
\begin{equation*}
\varphi (\widehat{x}+\eta^j ,\widehat{t})-\varphi (\widehat{x},\widehat{t}%
)-\langle D\varphi (\widehat{x},\widehat{t}),\eta^j \rangle =\langle
D^{2}\varphi (\widehat{x},\widehat{t})\eta^j ,\eta^j \rangle +\mathcal{O}(|\eta^j|^{2}).
\end{equation*}
Here and in the following,
we have used the notation $\eta^j = \eta^j(\widehat{x}, \widehat{t}, z)$.
Hence
\begin{align*}\label{eq:three-integs}
\mathcal{J}_j(\widehat{x},\widehat{t},D\varphi\left(\widehat{x},\widehat{t}\right),\varphi \left( \cdot , \widehat{t}\right) )
&= \int_{\mathbb{R}^{l}\setminus \{0\}}  \Big(\varphi \left( \widehat{x}+\eta^j, \widehat{t} \right) - \varphi \left( \widehat{x}, \widehat{t}\right) -\chi_{\{|z|\leq1\}} \langle \eta^j, D\varphi\left(\widehat{x},\widehat{t}\right) \rangle\Big) \nu^j \left(dz\right) \notag\\
&=\int_{B_{1}\left( 0\right) \setminus \{0\}}\mathcal{O}(|\eta^j|^{2})\nu^j(dz)
+ \int_{B_{1}\left( 0\right) \setminus \{0\}}\Vert D^{2}\varphi (\widehat{x},\widehat{t})\Vert |\eta^j |^{2}\nu^j(dz)\notag\\
&\quad+ \int_{\mathbb{R}^{\ell }\setminus B_{1}\left( 0\right) } \Big( \varphi (\widehat{x}+\eta^j ,\widehat{t})-\varphi (\widehat{x},\widehat{t}) \Big)\nu^j (dz) =: I_1 + I_2 + I_3
\end{align*}
and to bound $\mathcal{J}_j$ it suffices to establish an appropriate upper bound on the integrals $I_1, I_2$ and $I_3$.
Assumption \eqref{assump1elelnloc} yields
\begin{equation*}
I_1 = \int_{B_{1}\left( 0\right) \setminus \{0\}}\mathcal{O}(|\eta^j|^{2})\nu^j(dz) \leq C \left(1 + |\widehat{x}|^2 \right) \int_{B_{1}\left( 0\right) \setminus \{0\}} |z|^{2} \nu^j(dz) \leq C \left(1 + |\widehat{x}|^2 \right),
\end{equation*}
while assumption \eqref{assump1elelnloc} and \eqref{eq:secondderivativebound} gives
\begin{align*}
I_2 &=\int_{B_{1}\left( 0\right) \setminus \{0\}}\Vert D^{2}\varphi (\widehat{x},\widehat{t})\Vert |\eta^j |^{2}\nu^j(dz)\notag\\
&\leq c\, \lambda\, e^{c(s-\widehat{t})} C \left( 1 + |\widehat{x}|^{p} + |y|^{p-2}\right)\int_{B_{1}\left( 0\right) \setminus \{0\}} |z|^{2}\nu^j(dz)\notag\\
&\leq c\, \lambda\, e^{c(s-\widehat{t})} C \left( 1 + |\widehat{x}|^{p} + |y|^{p-2}\right).
\end{align*}
To estimate $I_3$ we observe that
$$
\varphi (\widehat{x}+\eta^j ,\widehat{t})-\varphi (\widehat{x},\widehat{t})
\leq \psi _{j}^{i,y,s,\delta}(\widehat{x}+\eta^j ,\widehat{t}) - \psi _{j}^{i,y,s,\delta}(\widehat{x},\widehat{t})
$$
since $\psi _{j}^{i,y,s,\delta} - \varphi$ has a minimum at $(\widehat{x}, \widehat{t})$.
Therefore
\begin{align*}
I_3 &= \int_{\mathbb{R}^{\ell }\setminus B_{1}\left( 0\right) } \Big(\varphi (\widehat{x}+\eta^j ,\widehat{t})-\varphi (\widehat{x},\widehat{t})\Big)\nu^j (dz)
\leq \int_{\mathbb{R}^{\ell }\setminus B_{1}\left( 0\right) } \Big(\psi _{j}^{i,y,s,\delta}(\widehat{x}+\eta^j ,\widehat{t}) - \psi _{j}^{i,y,s,\delta}(\widehat{x},\widehat{t})\Big) \nu^j (dz)\notag\\\notag\\
&\leq c\, \lambda\, e^{c(s-\widehat{t})} \int_{\mathbb{R}^{\ell }\setminus B_{1}\left( 0\right) }
\Big(\left(1 + |y|^{p-2}\right)\left( |\widehat{x} + \eta^j -y|^{2} - |\widehat{x}-y|^{2}\right) + |\widehat{x} + \eta^j -y|^{p} - |\widehat{x}-y|^{p}\Big)
\nu^j (dz)\notag\\
&\quad + \int_{\mathbb{R}^{\ell }\setminus B_{1}\left( 0\right) } \Big(c_{ij}(\widehat{x} + \eta^j,\widehat{t}) - c_{ij}(\widehat{x},t)\Big) \nu^j (dz).
\end{align*}
Using \eqref{assump1elelnloc} and \eqref{ass:reg} we obtain
\begin{align*}
I_3&\leq c\, \lambda\, e^{c(s-\widehat{t})} C \left( 1 + |\widehat{x}|^{p} + |y|^{p} \right) \int_{\mathbb{R}^{\ell }\setminus B_{1}\left( 0\right) }\nu^j (dz) + c\, \lambda\, e^{c(s-\widehat{t})} \left(1 + |y|^{p-2}\right) \int_{\mathbb{R}^{\ell }\setminus B_{1}\left( 0\right) } |\eta^j|^2 \nu^j (dz)\notag\\
&\quad+ c\, \lambda\, e^{c(s-\widehat{t})} \int_{\mathbb{R}^{\ell }\setminus B_{1}\left( 0\right) } |\eta^j|^p \nu^j (dz) + C \int_{\mathbb{R}^{\ell }\setminus B_{1}\left( 0\right) } \Big(1 + |\widehat{x}|^p + |\eta^j|^p \Big) \nu^j (dz)\notag\\
&\leq c\, \lambda\, e^{c(s-\widehat{t})} C \left\{ \left( 1 + |\widehat{x}|^{p} + |y|^{p} \right)  
+  |y|^{p-2} \left(1 + |\widehat{x}|^2\right) 
+  \left(1 + |\widehat{x}|^p\right)\right\}
+ C \left(1 + |\widehat{x}|^p\right)\notag\\
&\leq c\, \lambda\, e^{c(s-\widehat{t})} C \left( 1 + |\widehat{x}|^{p} + |y|^{p} \right).
\end{align*}
Summing up bounds for $I_1$, $I_2$ and $I_3$ implies
\begin{align}\label{eq:J-bound}
- \mathcal{J}_j(\widehat{x},\widehat{t},D\varphi\left(\widehat{x},\widehat{t}\right),\varphi \left( \cdot , \widehat{t}\right) ) \geq - c\, \lambda\, e^{c(s-\widehat{t})} C \left( 1 + |\widehat{x}|^{p} + |y|^{p} \right),
\end{align}
where $C$ may depend only on $K, n$ and $p$.

Finally, assumption \eqref{assump1elel} yields $|f^i(\widehat{x},\widehat{t})| \leq C \left( 1 + |\widehat{x}|^{p} \right)$ and
plugging this inequality, \eqref{eq:dt-bound}, \eqref{eq:L-bound} and \eqref{eq:J-bound} into \eqref{eq:eqpart} gives us
\begin{align}\label{eq:final-bound}
&-\partial _{t}\varphi \left( \widehat{x},\widehat{t}\right) +\mathcal{F}_{j}\left( \widehat{x},\widehat{t},\psi _{j}^{i,y,s}\left( \widehat{x},\widehat{t}\right),D\varphi \left( \widehat{x},\widehat{t}\right) ,D^{2}\varphi \left(\widehat{x},\widehat{t}\right) ,\varphi \left( \cdot ,\widehat{t}\right) \right) \notag\\
&\geq c^2\,\lambda\, e^{c (s-\widehat{t})}  \frac{1}{2^{p+2}} \left\{ 1 + |\widehat{x}|^{p} + |y|^{p} \right\}\notag\\
&\quad - c\, \lambda\, e^{c(s-\widehat{t})} C \left( 1 + |\widehat{x}|^{p} + |y|^{p} \right)
- C\left( 1 + |\widehat{x}|^p \right) \geq 0,
\end{align}
where the last inequality is based on choosing $c$ large enough, our choice depending only on $K, n$ and $p$.
Hence, we conclude that $\psi ^{i,y,s}$ satisfies \eqref{eq:eqpart} and therefore
the equation of being a viscosity supersolution is fullfilled.

This completes step 3 in the proof of Lemma \ref{lemma:supersol} and in light of step 1 and step 2 the proof of the Lemma is complete. $\hfill\Box$\\

\begin{remark}\label{remark:swcosts}
In the setting of optimal switching problems, assumption
\eqref{ass:structur} is no restriction.
In particular, we may assume, without loss of generality, that
\begin{equation}
c_{ii_{2}}(x,t)+c_{i_{2}i_{3}}(x,t)+\dots c_{i_{l-1}j}(x,t)\geq c_{ij}(x,t),
\label{eq:structureswcosts}
\end{equation}
for all switching chains $(i,i_{2},\dots ,i_{l-1},j)\in \mathcal{A}_{ij}$
and any $i,j\in \mathcal{I}_{m}$.
Indeed, if \eqref{eq:structureswcosts}
does not hold, we can construct new switching costs $\widetilde{c}_{ij}$ by
\begin{equation*}
\widetilde{c}_{ij}(x,t)=\min_{(i_{1},\dots ,i_{l})\in \mathcal{A}_{ij}}\sum_{k=1}^{l-1}c_{i_{k}i_{k+1}}(x,t)
\end{equation*}
which we then consider in place of $c_{ij}$.
Since the regularity assumptions \eqref{ass:reg} on the original switching costs $c_{ij}$ only assume semi-concavity
and Lipschitz continuity,
the new switching costs $\widetilde{c}_{ij}$ will satisfy \eqref{ass:reg} by construction.
The same is true for the classical no-loop condition in \eqref{ass:noloop}.
Moreover, using $\widetilde{c}_{ij}(x,t)$ in place of $c_{ij}$ will not alter the cost structure of the
problem and hence the value function will remain unchanged.
In the setting of optimal switching, the switching cost $\widetilde{c}_{ij}(x,t)$ represent switching using the
\textquotedblleft cheapest" switching chain from state $i$ to $j$ and
assumption \eqref{eq:structureswcosts} means that it is always cheaper to switch
directly to a state than to go through some intermediate state. More
explicitly, this implies that at any time $t$ at most one switch is made.
Note that the switching costs $\widetilde{c}_{ij}$ can be
no more than Lipschitz continuous, regardless of the regularity of $c_{ij}$.
Hence, it is essential that the regularity assumptions of Lundstr\"{o}m,
Nystr{\"{o}}m, and Olofsson \cite{LNO14b} are relaxed in order
to assume \eqref{eq:structureswcosts} without loss of generality in the context of optimal switching.
\end{remark}

\noindent
{\bf Proof of Theorem \ref{thm:time} (H\"older continuity in time)}
To prove H\"{o}lder continuity in the time variable, let $u$ be a
viscosity solution to \eqref{problem} and fix arbitrary $\left( y,s\right)
\in \mathbb{R}^{n}\times \left[ 0,T\right] $ and $i\in \mathcal{I}_{m}$. We
now apply Lemma \ref{lemma:supersol} with $h_{i}(x)=u_{i}(x,s)$. Note that
since $u$ is a viscosity solution to \eqref{problem}, Corollary \ref%
{cor:spatialreg} asserts that $h_{i}$ satisfies \eqref{eq:gass2}. By the
comparison principle, we have
\begin{equation*}
u_{j}(x,t) \leq \psi _{j}^{i,y,s}(x,t)
\end{equation*}
for all $j\in \mathcal{I}_{m}$ and $(x,t)\in \mathbb{R}^{n}\times \lbrack 0,s]$.
In particular, setting $j=i$ and $x=y$ and using $c_{ii}=0$,
this reduces to
\begin{equation*}
u_{i}(y,t)-u_{i}(y,s)\leq  c^2 e^{c(s-t)} A \left(\lambda \left(s - t\right) + \frac{1}{\lambda} \right)
\end{equation*}
where $A = \left( 1 + |y|^{p}\right)$.
For $s$ and $t$ fixed, we set $\lambda =$ $(s-t)^{-1/2} $ and get
%
\begin{equation*}
u_{i}(y,t)-u_{i}(y,s)\leq 2 c^2 e^{c(s-t)} A (s-t)^{1/2}.
\end{equation*}
Without loss of generality we may assume that $|s-t|\leq 1$, and, therefore,
\begin{equation*}
u_{i}(y,t)-u_{i}(y,s)\leq C(1+|y|^{p})(s-t)^{1/2},
\end{equation*}
where $C$ may depend only on $K,n,p$ and $T$.
This proves the H\"{o}lder continuity in the time variable `from above'.

The lower bound on $u_{i}(y,t)-u_{i}(y,s)$ follows similarly by considering the subsolution
$\check{\psi}^{y,s}$ from Corollary \ref{lemma:subsol} as a barrier from below.
The proof of Theorem \ref{thm:time} is thus complete. $\hfill \Box$\\


%
\noindent
{\bf Proof of Theorem \ref{thm:existence} (Existence)}
We construct a viscosity solution to problem \eqref{problem} using Perron's
method. To this end, we define $u=(u_{1},\dots ,u_{m})$ as
\begin{equation*}
u_{i}\left( x,t\right) :=\inf \{u_{i}^{+}\left( x,t\right) : u_{i}^{+} = (u_{1}^{+},\dots
,u_{m}^{+})\text{ is a supersolution to \eqref{problem}}\}
\end{equation*}
Note that this construction is well defined given the explicit viscosity
sub- and supersolutions $\check{\psi}_{j}^{y,T}(x,t)$ and $\psi_{j}^{i,y,T}(x,t)$ constructed in Lemma \ref{lemma:supersol} and Corollary \ref{lemma:subsol}.
We now intend to prove that $u^{\ast }$ and $u_{\ast }$,
the upper- and lower semicontinuous envelopes of $u$, are, respectively, a
subsolution and a supersolution to \eqref{problem}. It then follows by the
comparison principle that $u^{\ast }\leq u_{\ast }$ and hence $u=u_{\ast
}=u^{\ast }$ is a viscosity solution to \eqref{problem}.

We first prove that $u^{\ast }$ satisfies the terminal condition of being a
subsolution. To this end, we make use of the explicit viscosity sub- and
supersolutions from Lemma \ref{lemma:supersol} again. Fix a component $i\in
\mathcal{I}_{m}$ and a point $y\in \mathbb{R}^{n}$. 
By construction of $u_{i}$ we have $u_{i}(x,t)\leq \psi _{i}^{j,y,T}(x,t)$ whenever $(x,t)\in
\mathbb{R}^{n}\times \lbrack 0,T]$. 
Moreover, since $\psi _{i}^{j,y,T}(x,t)$
is continuous, it follows that $u_{i}^{\ast }(x,t)\leq (\psi
_{i}^{j,y,T}(x,t))^{\ast }=\psi _{i}^{j,y,T}(x,t)$ for every $\varepsilon >0$
and $(x,t)\in \mathbb{R}^{n}\times \lbrack 0,T]$. In particular, setting $%
j=i $ and $\left( x,t\right) =\left( y,T\right) $, we deduce
\begin{equation*}
u_{i}^{\ast }(y,T)\leq \lim_{\lambda \rightarrow \infty }\psi
_{i}^{i,y,T}(y,T)=u_{i}\left( y,T\right) +c_{ii}\left( y,T\right) =g_{i}(y).
\end{equation*}
Since $i$ and $y$ are arbitrary in this argument, we conclude that $u^{\ast
} $ satisfies the terminal condition. 
We prove that $u_{\ast}$ satisfies
the terminal condition in a similar way using the comparison principle and $\check{\psi}_{j}^{y,T}(x,t)$ from
Corollary \ref{lemma:subsol} as a barrier from below.

After noticing that our switching costs are continuous, the sub- and
supersolution properties are shown in the same way as outlined by Biswas,
Jakobsen, and Karlsen \cite[pages 70--72]{BJK10}. We refer the interested
reader there and simply conclude that the proof of Theorem \ref{thm:existence} is complete. $\hfill\Box$



\begin{thebibliography}{LNO14b}
\bibitem[AF12]{AF12} Brahim El-Asri and Imade Fakhouri, \emph{Optimal
multi-modes switching with the switching cost not necessarily positive},
arXiv:1204.1683v1, (2012).

\bibitem[AH09]{AH09} Brahim El-Asri and Sa\"{\i}d Hamad\`{e}ne, \emph{The finite
horizon optimal multi-modes switching problem: The viscosity solution
approach} Applied Mathematics \& Optimization, \textbf{60} (2009), 213-235.

\bibitem[BI08]{BI08} Guy Barles and Cyril Imbert, \emph{Second-order
elliptic integro-differential equations: viscosity solutions theory revisited%
}, Annales de l'Institut Henri Poincar\'{e}, \textbf{25} (2008), 567-585.

\bibitem[BJ07]{10BJK} Guy Barles and Espen R. Jakobsen, \emph{Error bounds
for monotone approximation schemes for parabolic Hamilton-Jacobi-Bellman
equations}, Math. Comp. \textbf{76} (2007), 1861-1893.

\bibitem[BJK10]{BJK10} Imran H. Biswas, Espen R. Jakobsen, and Kenneth H.
Karlsen, \emph{Viscosity solutions for a system of integro-PDEs and
connections to optimal switching and control of jump-diffusion processes},
Applied Mathematics \& Optimization, \textbf{62} (2010), 47-80.

\bibitem[BS85]{BS85} Michael J. Brennan and Eduardo S. Schwartz, \emph{%
Evaluating natural resource investments}, Journal of Business, \textbf{58} (1985),
135-157.

\bibitem[CIL92]{CIL92} Michael G. Crandall, Hitoshi Ishii, and Pierre-Louis
Lions, \emph{User's guide to viscosity solutions of second order partial
differential equations}, Bulletin of the American Mathematical Society,
\textbf{27} (1992), 1-67.

\bibitem[DHP10]{DHP10} Boualem Djehiche, Sa\"{\i}d Hamad\`{e}ne, and Alexandre
Popier, \emph{A Finite Horizon Optimal Multiple Switching Problem}, SIAM
Journal on Control and Optimization, \textbf{48} (2010), 2751-2770.

\bibitem[DHMZ14]{DHMZ14}
Boualem Djehiche, Sa\"{\i}d Hamad\`{e}ne, Marie-Amelie Morlais, and Xuzhe Zhao,
\emph{On the equality of solutions of max-min and min-max systems of
variational inequalities with interconnected bilateral obstacles},
arXiv preprint arXiv:1408.4282 (2014).

\bibitem[HM12]{HM12} Sa\"{\i}d Hamad\`{e}ne and Marie-Am\'{e}lie Morlais, \emph{%
Viscosity Solutions of Systems of PDEs with Interconnected Obstacles and
Switching Problem} Applied Mathematics \& Optimization, \textbf{67}(2013),
163-196.

\bibitem[HM16]{HM16} Sa\"{\i}d Hamad\`{e}ne and Marie-Am\'{e}lie Morlais,
\emph{Viscosity solutions for second order integro-differential equations without monotonicity condition:
the probabilistic approach},
Stochastics \textbf{88} (2016): 632-649.

\bibitem[HZ15]{HZ15} Sa\"{\i}d Hamad\`{e}ne and Xuzhe Zhao \emph{Systems of
integro-PDE with interconnected obstacles and multi-modes switching problem
driven by L\'{e}vy process}, Nonlinear differential equations and
Applications, \textbf{22} (2015), 1607-1660.

\bibitem[HL13]{HL13}
F. Reese Harvey, and H. Blaine Lawson, 
\emph{Notes on the differentiation of quasi-convex functions},
arXiv preprint arXiv:1309.1772 (2013).

\bibitem[HT07]{HT07} Ying Hu and Shanjian Tang, \emph{Multi-dimensional BSDE
with Oblique Reflection and Optimal Switching}, Probability Theory and
Related Fields, \textbf{147} (2010), 89-121.

\bibitem[IK91]{IK91} Hitoshi Ishii and Shigeaki Koike, \emph{Viscosity
solutions of a system of nonlinear second-order elliptic PDEs arising in
switching games}, Funkcialaj Ekvacioj, \textbf{34} (1991), 143-155.

\bibitem[IS04]{IS04} Hitoshi Ishii, Moto-Hiko Sato, \emph{Nonlinear bolique
derivative problems for singular degenerate parabolic equations on a general
domain}, Nonlinear analysis, \textbf{57} (2004), 1077-1098.

\bibitem[JK06]{JK06} Espen R. Jakobsen and Kenneth H. Karlsen, \emph{A Maximum
Principle for Semicontinuous Functions Applicable to Integro-Partial
Differential Equations}, Nonlinear differential equations and applications,
\textbf{13} (2006), 137-165.
ƒ
\bibitem[K16]{K16}
Tomasz Klimsiak,
\emph{Obstacle problem for evolution equations involving measure data and operator corresponding to semi-Dirichlet form}, arXiv preprint arXiv:1612.07274 (2016).

\bibitem[K16b]{K16b}
Tomasz Klimsiak,
\emph{Systems of quasi-variational inequalities related to the switching problem},
arXiv preprint arXiv:1609.02221 (2016).

\bibitem[Kr05]{34BJK} Nikolai V. Krylov, \emph{The rate of convergence of
finite-difference approximations for Bellman equations with Lipschitz
coefficients}, Appl. Math. Optim., \textbf{52} (2005) 365-399.

\bibitem[LNO14]{LNO14} Niklas L.P. Lundstr\"{o}m, Kaj Nystr{\"{o}}m and Marcus
Olofsson, \emph{Systems of variational inequlities for nonlocal operators
related to optimal switching problems: existence and uniqueness},
Manuscripta Mathematica, \textbf{145} (2014), 407-432.

\bibitem[LNO14b]{LNO14b} Niklas L.P. Lundstr\"{o}m, Kaj Nystr{\"{o}}m, and
Marcus Olofsson, \emph{Systems of variational inequalities in the context of
Optimal Switching Problems and Operators of Kolmogorov Type}, Annali di
Mathematica Pura ed Applicata, \textbf{193} (2014), 1213-1247.

\bibitem[L\"O15]{L�1�7}
Niklas L.P. Lundstr\"{o}m and Thomas \"Onskog,
\emph{Stochastic and partial differential equations on non-smooth time-dependent domains},
arXiv preprint arXiv:1503.05433 (2015).

\bibitem[M14]{M14}
Randall Martyr,
\emph{Finite-horizon optimal multiple switching with signed switching costs},
arXiv preprint arXiv:1411.3971 (2014).


\end{thebibliography}
\end{document}